\documentclass[11 pt]{amsart}

\usepackage{ amsfonts,amscd,amssymb,amsthm}
\usepackage{amsfonts}
\usepackage{amssymb}
\usepackage{amsmath}
\usepackage{amsxtra}
\usepackage{amscd}
\usepackage{amsthm}
\usepackage{eucal}
\usepackage{graphicx}

\def\Zz{{\mathbb Z}}

\title{On algebraic closure in pseudofinite fields}
\author{\"{O}zlem Beyarslan, Ehud Hrushovski}
\date{\today}
\newtheorem{theorem}{Theorem}
\newtheorem{definition}[theorem]{Definition}
\newtheorem{proposition}[theorem]{Proposition}
\newtheorem{lemma}[theorem]{Lemma}
\newtheorem{corollary}[theorem]{Corollary}

\newtheorem{Remark}[theorem]{Remark}
\newtheorem{problem}[theorem]{Problem}
\newtheorem{question}[theorem]{Question}
\newtheorem{example}[theorem]{Example}
\newtheorem{notation}[theorem]{Notation}

\def\nn {\mathbb{N}}
\def\Nn{{\mathbb N}}
\def\zz {\mathbb{Z}}
\def\qq {\mathbb{Q}}

\newfont{\mym}{msbm10}

\def\1{^{\circ}}
\def\2{_{\circ}}

\def\claima{\noindent {\bf Claim 1:} $\,$}
\def\claimb{\noindent {\bf Claim 2:} $\,$}
\def\proof{\noindent {\bf Proof:} $\,$}

\def\-#1{\overline{#1}}
\def\a{\alpha}

\def\w{\omega}
\def\si{\sigma}

\def\G{\Gamma}

\def\ontor{\rlap{$\to$}\kern.05pt\to}

\def\union{\cup}
\def\meet{\cap}

\def\qed{$\Box$}

\def\ch{\mathop{\rm char}\nolimits}

\def\acl{\mathop{\rm acl}\nolimits}

\def\Sym{\mathop{\rm Sym}\nolimits}

\def\Aut{\mathop{\rm Aut}\nolimits}

\def\Ker{\mathop{\rm Ker}\nolimits}
\def\Gal{\mathop{\rm Gal}\nolimits}

\def\mod{\mathop{\rm mod}\nolimits}

\def\dcl{\mathop{\rm dcl}\nolimits}
\def\acl{\mathop{\rm acl}\nolimits}
\def\Psf{\mathop{\rm Psf}\nolimits}

\def\mod{\mathop{\rm mod}\nolimits}
\def\qed{\hfill$\Box$}
\def\proof{\noindent {\bf Proof:} $\,$}

\def\inv{ ^{-1}}
\thanks{The second author was supported by ISF 1048/07}

\begin{document}

\maketitle
\date{today}

\begin{abstract}
We study the  automorphism group of the algebraic closure of a
substructure $A$ of a pseudo-finite field $F$.  We show that the
behavior of this group, even when $A$ is large, depends essentially
on the roots of unity in $F$.  For almost all completions of the
theory of pseudofinite fields, we show that over $A$, algebraic
closure agrees with definable closure, as soon as $A$ contains the
relative algebraic closure of the prime field.
\end{abstract}

\section{Introduction}

A pseudofinite field is an infinite model of the theory of finite fields.
By Ax \cite{Ax}, we know that a field $F$ is pseudofinite if and only if it is 1) perfect, 2) PAC and 3) has a unique (and so
necessarily Galois and cyclic) extension in the algebraic closure $F^a$ of $F$ of degree $n$ for every $n\in \nn\setminus\{0\}$.
See   \cite{jardenfried} for the PAC property; it  will play almost no role in this paper.

We are interested in definable and algebraic closure in $F$, over a substructure $A$ containing an elementary submodel $M$.   Surprisingly, the answer depends
intimately on embeddings of number fields, or finite fields, into $M$.  In characteristic  zero for instance, we show that model-theoretic Galois groups over $M$
have odd order, unless $M$ contains all even-order roots of unity.

Real closed fields provide a geometrically comprehensible way of  symmetry-breaking in algebraic geometry; a Galois cover of an algebraic variety splits into semi-algebraic sections.   Our results imply that ``almost all" pseudo-finite fields
give an alternative geometric approach to such a splitting:   the Galois cover splits into definable sections, or in model theoretic terms, definable and algebraic closures coincide.  See \S 5, Corollary 8.

These symmetry-breaking results  are in fact valid for quasi-finite fields in the
sense of \cite{Serre}, p. 188, i.e. perfect fields with absolute
Galois group $\hat{\Zz}$, the profinite completion of $\Zz$.  We use
pseudo-finiteness only in order to demonstrate the converse, that if
all $p^n$'th roots of unity are contained in  $F$, then Galois groups of order divisible by $p$ do occur as $Aut(B/A)$
with $M \leq A \leq B \leq N$, $M \prec  N$.

In \S \ref{tourn} ,we describe (in characteristic prime to $p$) an explicit structure incompatible with symmetries of order $p$.

We thank the referee for a close reading and two very useful reports.

\section{Quasi-finite Fields}

We write $\leq$ for the substructure relation; in particular, for fields $A,B$, $A\leq B$ means that $A$ is a subfield of $B$.  The algebraic closure of any field $F$ will be denoted $F^a$.

By definition, a quasi-finite field $F$ has a unique
 extension in $F^a$ of degree $n$ for every $n\in \nn\setminus\{0\}$. Let $F_n$ denote the unique extension of $F$ in $F^a$ of degree $n$. This extension is easily seen to be interpretable in $F$ using
parameters from $F$. Indeed, as $F$ is perfect, $F_n=F(\a)$ for some
$\a \in F_n$. Let $X^n+a_1X^{n-1}+\cdots +a_n$ be the minimal
polynomial of $\a$ over $F$. Then $F_n$, which is an $n$-dimensional
vector space over $F$ with basis $\{1,\a,\a^2,\ldots,\a^{n-1}\}$, is
definably isomorphic to $F^n$ via this basis as a vector space. Also
any linear endomorphism of $F_n$ translates into a definable (with
parameters) linear endomorphism of $F^n$ (coded by an $n\times n$
matrix over $F$). In particular, the $\a$-multiplication in $F_n$,
the multiplicative structure of $F_n$ and the action of
$\Gal(F_n/F)$ on $F_n$ can all be definably (with parameters) coded
in $F^n$.

Note that to interpret $F_n$ in $F$ we only need $a_1,\,\ldots, a_n$ as parameters, but to interpret the action of an element $\tau$ of
$\Gal(F_n/F)$ in $F$, apart from these $n$ parameters, we also need $b_0, \ldots, b_{n-1} \in F$ where $\tau(\a) = b_0 + b_1\a + \cdots +
b_{n-1}\a^{n-1}$, which makes up a total of $2n$ parameters. Note also that any other choice of the parameters $a_1,\,\ldots, a_n$ for which the
polynomial $X^n+a_1X^{n-1}+\ldots +a_n$ is irreducible gives rise to an isomorphic structure $F_n$; on the other hand, different choices of the
parameters $b_0,\,\ldots, b_{n-1}$ may define different field automorphisms.

\begin{lemma} \label{commute}Let $F$ be a quasi-finite field and $\si$ be a topological generator of $\Aut(F^a/F)$, let $M$ be an elementary submodel of $F$. Let $\mu$ be in $\Aut(F/M)$. Then any extension of $\mu$ to
$F^a$ commutes with $\sigma$.
\end{lemma}
\proof It is enough to show that $\sigma$ and $\mu$ commute on $F_n$ where $F_n$ is the unique extension of $F$ of degree $n$. Since $M$ is an elementary submodel of $F$,
$F_n=F(\alpha)$ where $\alpha$ is a root of an irreducible polynomial of degree $n$ with coefficients in $M$.  We will show that $\sigma \mu (\alpha)= \mu \sigma (\alpha)$. We know that $\sigma(\a)=b_0+b_1\a+\ldots +b_{n-1}\a^{n-1}$ for some $b_0,\ldots b_{n-1}$ in $M$. Then $\mu(\sigma(\a))=b_0+b_1\mu(\a)+\ldots +b_{n-1}\mu(\a)^{n-1}$. On the other hand since $\mu$ is a field automorphism fixing the minimal polynomial of $\a$, $\mu(\alpha)=\sigma^r(\alpha)$ for some $r<n$, hence $\sigma(\mu(\a))=\sigma(\sigma^r(\alpha))=\sigma^r(\sigma(\alpha))=\sigma^r(b_0+b_1\a+\ldots +b_{n-1}\a^{n-1})=b_0+b_1\sigma^r(\a)+\ldots +b_{n-1}\sigma^r(\a)^{n-1}=b_0+b_1\mu(\a)+\ldots +b_{n-1}\mu(\a)^{n-1}=\mu(\si(\alpha))$.

We include here also a lemma of homological flavor, that will be essential in the main theorem.

\begin{lemma} \label{essential}Let $D$ be an abelian group with three commuting endomorphisms $P,S,T$. Let $\Omega=\cup_n \ker{P^n}$. Assume:
\begin{enumerate}
\item $P$ is surjective.
\item $T|\Omega =0$
\item $\Omega\cap \ker(S)\subseteq  \ker(P^r)$ for some $r \in \Nn$.

Then: if $a\in \ker(S)$ and $P(a)\in \ker(T)$, then $a\in \ker(T)$.
\end{enumerate}
\end{lemma}
\proof Let $P(a)=b$ and $C=\{x\in D :   \, P^n(x)=b\ \hbox{for some
} \, n > 0\}$. Since $S(b)=T(b)=0$ we have $S(C)$, $T(C)\subseteq
\Omega$. By (ii) $TS | _C=0$ i.e. $T(C)\subseteq  \ker(S)$.  By (iii), $T(C) \subseteq
\ker(P^r)$. But $P(C)=C$   by (i) and the definition of $C$, so
$PT(C)=T(C)$, hence $P^rT(C)=T(C)$; so $T(C)=0$.
\endproof

\subsection{Notation:  maximal $p$-extensions, roots of unity}

Let $k$ be a prime field,  $p$ a prime.  If  $p\neq \ch(k)$, we let $\mu_{p^n}$ denote the multiplicative subgroup of
$k^a$ of $p^n$-th roots of unity; and $\mu_{p^\infty}=\bigcup\limits_{n<\omega} \mu_{p^n}$.  Also write $\Omega_p = \mu_{p^\infty}$.

On the other hand, if $p= \ch(k)$, we let $\Omega_p$ be the  maximal $p$-extension
of the prime field; in other words $\Omega_p = \union_n \Ker(P^n)$, where $P(x)=x^p-x$
is the Artin-Schreier operator.

   In
 either case it is clear that any field $M$ either contains $\Omega_p$, or intersects
 $\Omega_p$ in a finite group.

\section{Geometric Representation}

In this section we will state and prove our main theorem on automorphism groups of pesudofinite fields.  We begin with the main definitions.
Let $G$ be a pro-finite group.

\begin{definition} We say that the group $G$ is \emph{geometrically represented} in the theory $T$ if there exists $M_0\prec M \vDash T$ and $M_0 \leq A \leq B \leq M$,   such that   $B \subseteq \acl(A)$ and $\Aut(B/A)  \cong G$.  We say that a prime number $p$ is geometrically represented in the theory $T$ if $p$ divides the order of  some finite group $G$ geometrically represented in $T$.
\end{definition}

In this definition, $A,B$ are substructures of $M$ containing $M_0$.
$\Aut(B/A)$ must be interpreted as the set of permutations of $B$
over $A$ preserving the truth value of all formulas (computed in
$M$.)  If one takes $M$ to be saturated and of greater cardinality than $B$, this can also be described as  the
set of   permutations of $B$ fixing $A$ that extend to   automorphisms of $M$.
Compare \cite{H}.

For a theory of perfect fields,  by Galois theory, $p$ is geometrically represented in $T$ iff the group $\Zz/p\Zz$ is  geometrically represented in $T$.    As the referee pointed out, for general theories this may not hold, and one may prefer a definition allowing $A,B$ to consist of imaginary elements.  Theories of pseudo-finite fields admit elimination of finite imaginaries
over a model $M$ (cf. \cite{PAC}), so for such theories the two options are the same.  For simplicity and as we are only concerned with fields, we will use the definition above.

\begin{Remark}  \label{product1}    If a finite group $G$ is geometrically represented in the complete theory $T$
over a model $M_0$, then $G$ is also geometrically represented over any elementary  extension $M$ of $M_0$.  {\rm  Indeed we may assume
$G=\Aut(B/A)$ where $B=M_0(ab),A=M_0(a)$.    For any enumeration $m$ of $M$, $tp(m/M_0)$ has an extension to $M_0(a,b)$ 
which is finitely satisfiable in $M_0$; in this situation one says that $tp(ab/M)$ is a heir of $tp(ab/M_0)$.  So me may take 
  $tp(ab/M)$ to have this property.  It follows that the multiplicity of $b$ over $M(a)$ cannot be smaller than that of $b$ over $M_0(a)$.
So $tp(b/M_0(a)) \vdash tp(b/M(a))$.  Hence $G=\Aut(M(a,b)/ M(a))$.}   \end{Remark}

\begin{theorem}\label{aut} Let $F$ be a quasifinite field, $p$ a prime; if $p \neq \ch(F)$, assume  $F$ contains a primitive $p$'th root of unity.  Assume
$p$ is geometrically represented in $Th(F)$.  Then    $F$  contains the group  $\Omega_p$, $p^n$'th roots of unity if $p \neq \ch(F)$, maximal $p$-extension of the prime field if  $p = \ch(F)$.
\end{theorem}

\proof   By assumption there exist $M \prec N \models Th(F)$,
$M \leq A   \leq B$,   such that  $p$ divides
$|\Aut(B/A)|$.       Replacing $A$ by $Fix(\tau_B)$,
where $\tau_B$ is some element in $\Aut(B/A)$ of order $p$, we
 may assume $B/A$ is a Galois extension
of order $p$, generated by $\tau_B$.  We may take $N$ to be $|B|^+$-saturated.
Let $\tau_N$ be an extension of $\tau_B$ to $\Aut(N)$.
 Since $M$ is an elementary submodel of $N$, $M$ is relatively algebraically closed in $N$, $N$ and $M^a$ are
linearly disjoint over $M$; hence we may extend $\tau_N$ to a field automorphism
 ${\tau}$ of $N^a$, in such a way that ${\tau}$ fixes $M^a$.

 Since $F$ is quasi-finite, we have $\Aut(M^a/ M) =\Aut(N^a/N) \cong \widehat{\Zz}$;
 there exists an automorphism $\si$ of $N^a/N$ generating the Galois group $\Aut(N^a/N)$; the restriction $\si |M^a$ generates   $\Aut(M^a/M)$.  By Lemma \ref{commute},
  ${\tau}$ commutes with $\si$.

 Let $G$ be the multiplicative group if $p$ is not the characteristic of $F$, and let $G$ be the additive
 group otherwise.  Let $(D,+)= G(N^a)$, written additively.    $End(D)$ is also written additively.   Let $S = \si-Id$ and $T=\tau - Id$;
 these are commuting endomorphisms of $D$.  We define an additional endomorphism $P$ commuting with both, and an element $a \in A$,  according to cases:

 \begin{itemize}
 \item  If $p=\ch(F)$, let $P(x)$ be the Artin-Schreier operator $Fr(x)-x$, where $Fr$ is the Frobenius $p$'th power map on $G$.
 In this case, by Artin-Schreier theory, $B=A(b)$ for some $b$ with $P(b) \in A$.
  \item If $p \neq \ch(F)$ and $M$ contains the $p$'th roots of $1$, let $P=p \in End(D)$, $P(x)=px$. By Kummer theory, $B=A(b)$
  for some $b$ with $P(b)\in A$.   
 \end{itemize}

   Define $ \Omega$ be as in Lemma \ref{essential}; so $\Omega=\Omega_p$.    As $\si, {\tau}$
 are field automorphisms, they commute with $P$; so $P,S,T$ commute.

  (i) $P$ is
surjective on $N^a$ by algebraic closedness.

 (ii) $T|\Omega_p=0$ since $\tau$ fixes $M^{a}$.

 (iii) Suppose for contradiction that $\Omega_p$ is {\em not} contained
in $M$.  Then $M \meet \Omega_p$ is a finite subgroup of $\Omega_p$,
and for some $r \in \Nn$,  $P^r$ vanishes on $M \meet \Omega_p = \ker(S) \meet \Omega_p$.

By Lemma \ref{essential},  $T(b)=0$, i.e. $\tau(b)=b$; so $\tau |B=Id_B$.   This contradicts the choice of $\tau$.
    Thus $\Omega_p$ is contained in $M$, and so in $F$.
\qed
\endproof
\vskip 0.5cm

 \begin{corollary} \label{aut-c}  Let $F$ be a quasifinite field, $p$ a prime; assume $F[\mu]$ does not contain $\Omega_p$, where $\mu$ is a primitive $p$'th root of $1$.
 Then $p$ is not geometrically represented in $Th(F)$.
 \end{corollary}

 \proof Let $F'=F[\mu]$.  Since $[F':F] $ divides $p-1$, it is clear that $p$ remains  geometrically represented in $Th(F')$.
  By Theorem \ref{aut}, $F'$ contains $\Omega_p$.
   \qed
\endproof
\vskip 0.5cm
   We prove a converse   to Theorem \ref{aut} when $\ch(F) \neq p$, and $F$ is pseudo-finite.

\begin{theorem}\label{aut-conv} Let $p$ be a prime, $F$ a pseudofinite field not of characteristic $p$.  Assume
$F$  contains $\mu_{p^\infty}$.   Then
  $p$ is geometrically represented in $Th(F)$.
\end{theorem}
\proof

As $Th(F)$ is pseudo-finite, it is the restriction to $Fix(\si)$
of a completion $T$ of the theory ACFA of algebraically closed fields with an automorphism $\si$.
  We refer to \cite{CH} for basic facts about ACFA.  In particular,
   If $A$ is a substructure of a model of $T$ and $\acl(A)=A$,  then any automorphism $\tau$ of $(A,\si)$ is elementary; so $\tau$ restricts to an automorphism of $Fix(\si)$, elementary
in the sense of $Th(F)$.

Let $K=Fix(\si)$ where $(M,\si) \models T$.  (One may choose $M$ countable, if desired.)
  
Let $N$ be the field of generalized power series in $x$ with $\qq$-exponents
with coefficients in $M$. By \cite{Hahn} this  is an algebraically
closed field, see \cite{kedlaya}. Extend $\si$ to $N$ by mapping
$\sum \alpha_i x^i$ to $\sum \si(\alpha_i) x^i$.  Then $(N,\si)$
embeds into an elementary extension of $(M,\si)$.

Let $\{\w_i\}_{i<\omega}$ be a coherent system of the $p^i$-th roots of unity in $M$, \emph{i.e.} $\w_0=1$ and $\w_{i+1}^p=\w_i$ for $i\geq 0$.
Define $\tau$ to be an automorphism of $N$ fixing  $M$, and acting naturally on generalized power series, via:
 $$\tau:
x^{1/p^i}\mapsto \omega_i x^{1/p^i}$$ and
$$\tau: x^{1/n}\mapsto x^{1/n} \hbox{ for }  p\nmid n.$$

Note that, for $\si(x^{1/p^i})= x^{1/p^i}$ we have that  $$\si(\tau(x^{1/p^i}))=\si(\w_i x^{1/p^i})=\w_i   x^{1/p^i}$$
$$\tau(\si(x^{1/p^i}))=\tau( x^{1/p^i})=  \w_i x^{1/p^i},$$
As $\si,\tau$ also commute on $x^{1/n}$ for $(n,p)=1$, and on $M$, it is clear that   $\si$ commutes with $\tau$ on $N$.

  Now $\tau$ fixes $K(x)$ but not $K(x^{1/p})$; and $\tau^p$ fixes $K(x^{1/p})$.
 So the group of $N$- or $K$-elementary automorphisms of $K(x^{1/p})$ over $K(x)$ includes (and hence equals) $\langle\tau | K(x^{1/p})\rangle \cong \Zz/p$.

\section{ Automorphism Group and Tournaments} \label{tourn}

Here we give a different proof of  Theorem \ref{aut}, by constructing a structure that can have no automorphisms of order $p$.

Let $p$ be a prime.  By a {\em $p$-tournament} we mean a $p$-place relation $R$, such that
for any $p$-tuple of distinct elements $x_1,\ldots,x_p$,
$$R(x_{\tau(1)},\ldots,x_{\tau(p)}) \hbox{ holds for exactly one element  }  \tau \in \langle( 12\ldots p)\rangle$$
where $(12\ldots p)$ denotes the cyclic permutation of order $p$ over the $p$ element set $\{1,\ldots,p\}$,
and $ \langle( 12\ldots p)\rangle$ is the subgroup of $\Sym(p)$ generated by this permutation, isomorphic to $ \mathbb{Z}/p\mathbb{Z}.$

A $p$-tournament clearly has no automorphism of order $p$, or even an automorphism $\si$ with a $p$-cycle
$a_0,\ldots,a_{p-1}$ with $\si(a_i)=a_{i+1} (\mod p)$.   Thus $p$ is {\em not} geometrically represented in $T$
if $T$ is 1-sorted and admits a $p$-tournament structure on the main sort.  In fact no Galois group of $T$ can
have order divisible by $p$, whether or not the base contains an elementary submodel.

\def\m{\setminus}

\begin{proposition}  \label{p-tournament} Let $p$ be a prime, and  $F$ a field of characteristic $\neq p$,
containing the group $\mu_p$ of $p$'th roots of unity.   Let $\omega
\in \mu_p \m \{1\}$.  Let $S$ be a set of representatives for the
cosets of $\mu_p$ in $F^*=F\m \{ 0\}$.
  Then in the structure $(F,+,\cdot,\omega,S)$ there exists a  definable
$p$-tournament on $F$.  \end{proposition}

\begin{remark}  \rm
When $F$ is pseudo-finite, and $\omega \in F$, we have
$[F^*:(F^*)^p]=p$ by a counting argument.  The same conclusion holds
when $F$ is quasi-finite, using Galois cohomology: the cohomology
exact sequence associated with the short exact sequence  $$1 \to
\mu_p \to (F^{a})^*\to_{x \mapsto x^p} (F^{a})^* \to 1$$ gives,
using Hilbert 90,  $$F^* \to_{x \mapsto x^p} F^* \to
Hom(\widehat{\Zz}, \mu_p) \to H^1(Gal, (F^{a})^*)=0.$$ We refer to
\cite{Tate} for the basics of Galois cohomology.

 Assume $F$ contains
 a primitive $p^n$-th root of unity  $\zeta$, but not any $p$'th root of $\zeta$.    Then
 $F^*$ is the direct sum of $\mu_{p^n}$ and $(F^*)^{p^n}$.  Let $S_0$ be a set of representatives
 for $\mu_{p^n} / \mu_p$.  Then $S_0 (F^*)^{p^n}$ is a set of representatives for $(F^*)/ \mu_p$.
 Hence, using the Proposition, there exists a $p$-tournament definable in the field $F$ using
 $\mu_{p^n}$ as parameters.  This gives another proof of Theorem \ref{aut}.
\end{remark}

Before proving Proposition \ref{p-tournament}, we illustrate it with  the   case $p=2$.
  Assume $F$ does not contain $\sqrt{-1}$.
A \emph{tournament} on a set $X$ is an irreflexive binary relation $R\subset X\times X$ such that for every $x\neq y \in X$ exactly one of
$R(x,y)$ and $R(y,x)$ holds. A pseudofinite field $F$ not containing $\sqrt{-1}$ interprets a tournament by the formula: $$(\exists z) (z^2 =
x-y).$$ The automorphism group of any field interpreting a 0-definable tournament can not have any involutions.

We can still define a tournament in a pseudofinite field $F$ which  contains all the $2^n$-th roots of unity but not all the $(2^{n+1})$-st
roots of unity.

For every $m\in \nn$ we denote the set of $2^m$-th roots of unity by $\mu_{2^m}$. Let  $S\subset \mu_{2^n}$, such that $S\meet -S=\emptyset$ and
$S\cup -S=\mu_{2^n}$. Define a relation $R$ on $F\times F$ as follows:
$$R(x,y) \hbox{ if and only if } x-y  \hbox{ is in } \bigcup_{c\in S}c F^{2^n}.$$

Then this defines a tournament in $F$. That is, for every $x,y \in F, \, x\neq y$, exactly one of $R(x,y)$ and $R(y,x)$ holds. Suppose $\neg
R(x,y)$ then $(x-y)\not\in \bigcup_{c\in S} c F^{2^n}$ then $(x-y)$ is in $\bigcup_{c\in -S} c F^{2^n}$. Therefore $$-(x-y)=(y-x) \in
\bigcup_{c\in S} c F^{2^n}$$ hence $R(y,x)$. Also, at most one of $R(x,y)$ and $R(y,x)$ hold since $$F^{*}=\bigsqcup_{c\in
\mu_{2^n}}c{F^{*}}^{2^n},$$ that is, $\mu_{2^n}$ is a set of representatives for the cosets of the subgroup ${F^{*}}^{2^n}$ of
multiplicative part $F^{*}$ of $F$.

Now we will generalize the construction of the above tournament relation from binary to $p$-ary.

\proof  (of Proposition \ref{p-tournament})

Define a $p$-ary relation $R_\omega$ on $F$ as follows:
$$R_\omega(x_1,x_2,\ldots,x_p) \hbox{ if and only if } x_1+\omega x_2+
\ldots+\omega^{p-1}x_{p} \in S $$
\claima  Assume $x_1+\omega x_2+
\ldots+\omega^{p-1}x_{p} \neq 0$.  Then
$$R_\omega(x_{\tau(1)},\ldots,x_{\tau(p)}) \hbox{ holds for exactly one element in }  \langle( 12\ldots p)\rangle\simeq \mathbb{Z}/p\mathbb{Z}.$$

Indeed let $\pi \in  \langle( 12\ldots p)\rangle$ and $k=\pi(1)$ (so $k$  determines the element $\pi$). Then we have:
$$x_{\pi(1)}+\omega x_{\pi(2)}+\ldots+\omega^{p-1}x_{\pi(p)}=\w^{-(k-1)}(x_1+\omega x_2+\ldots+\w^{p-1}x_{p})$$
Since $S$ is a set of representatives for $F^*/\mu_p$, and
$a:=x_1+\omega x_2+\ldots+\w^{p-1}x_{p} \in F^*$, it is clear that
$\w^{-(k-1)} a \in S$ for a unique value of $k$ modulo $p$.

Thus $R_\omega$ is almost a $p$-tournament, but we need to deal with certain linearly dependent $p$-tuples.

\claimb   Assume $x_1 + \omega^i x_2 + \cdots + \w^{i(p-1)} x_p = 0$ for all $i=1,\ldots,p-1$.  Then
$x_1=\cdots =x_p$.

This is because the Vandermonde matrix with rows $(1,\omega,\ldots,\omega^{p-1})$,
$(1,\omega^2,\cdots,\omega^{2(p-1)})$, $\ldots$, $(1,\omega^{p-1},\ldots,\omega^{(p-1)(p-1)})$
has rank $p-1$.  So the kernel of this matrix is a vector space of dimension $1$.  But
$(1,\ldots,1)$ is clearly in the kernel; hence the kernel consists of scalar multiples of this vector.

Since we are only concerned with $p$-tuples of distinct elements, for each such $p$-tuple
$x=(x_1,\ldots,x_p)$ there exists a smallest $i \in \{1,\ldots,p-1\}$ such that $x_1+ \omega^i x_2 + \ldots \neq 0$.
Write $i=i(x)$, and define $R(x_1,\ldots,x_p) $ to hold iff $R_{\omega^{i(x)}}(x_1,\ldots,x_p)$ holds.
It is then clear that $R$ is a $p$-tournament.\qed

\section{Model Theoretic Consequences}

 Let $T_{\Psf}$ be the theory of pseudo finite fields. Let $K= \mathbb{Q}$ or $K=\mathbb{F}_p$.  By Ax's theorem \cite{Ax} (cf. also  \cite{jardenfried}, Chapter 20)  there is a one to one correspondence between the conjugacy classes of $\Aut(K^a/K)$  and the set of completions of the theory $T_{\Psf}$ of characteristic $=\ch(K)$.  Namely, note that if $M$ is a model of $T$,   $K^a \meet M$ is determined by $T$ up to isomorphism; call it $K^a_T$.  Then
 $\si$ corresponds to $T$ iff $Fix(\si) \cong K^a_T$.

  The absolute Galois group $\Gamma=Gal(K^a/K)$ is a compact topological group with a unique normalized  left invariant Haar measure $\mu_\Gamma$.  Let $\Pi$ be the set of conjugacy
  classes of $\Gamma$, and let $\pi: \Gamma \to \Pi$ be the quotient map.
  $\mu_\Gamma$ induces a measure $\mu$ on $\Pi$, namely $\mu(U) = \mu_\Gamma(\pi\inv(U))$.     Using the  1-1 correspondence above, we identify $\Pi$ with the the set of completions $\mathcal{C}$ of the theory of pseudofinite fields of characteristic =$\ch(K)$.  We obtain a measure on $\mathcal{C}$.   By a theorem of Jarden (cf. Theorem 20.5.1 of \cite{jardenfried}),
  for almost all $\si \in \G$, $K^a_T \models T_{\Psf}$.

If $A \subset N \models T$, let  $\dcl(A)$ denote the definable closure of $A$ in $N$.    We also write $M(a)$ for $\dcl(M \union \{a\})$; if $A,B$ are definably closed subsets of $N$, we will just write $AB$
for $\dcl(AB)$.
 $\acl$ denotes algebraic closure.   Thus $b \in \acl(A)$ iff there exists a formula $\phi(x)$ with parameters in $A$ such that $N \models \phi(b)$ and
 $\phi(N)$ is finite.  The smallest possible size $|\phi(N)|$ is called the multiplicity of $b$ over $A$.

 \begin{corollary} For almost all $T$ in $\mathcal{C}$,  we have $\acl=\dcl$ over
 $K^a_T$.
 \end{corollary}
 \proof For each prime $p\neq \ch (K)$ the set $\{\sigma \in \Aut (K^a/K) : \sigma^{p-1} \hbox{ fixes } \, \mu_{p^{\infty}} \}$ has measure 0. So $\bigcup_{p\neq \ch (K)} \{\sigma \in \Aut (K^a/K) : \sigma^{p-1} \hbox{ fixes } \, \mu_{p^{\infty}} \}$ has measure 0. If $\ch (K)=p_0$, $\{\sigma \in \Aut (K^a/K) : \sigma \hbox{ fixes the maximal }  p_0 \hbox { extension } L_{p_0} \}$ has measure 0.   Hence by Corollary  \ref{aut-c}, for almost all $T\in \mathcal{C}$, any group which is geometrically represented in $T$ is trivial;  hence $acl=dcl$ over  $K^a_T$.

Clearly, the same is true for Baire category in place of measure.

\begin{Remark} \rm
  While $\dcl=\acl$ is a restricted form of
Skolemization, the theories of pseudo-finite fields are not
Skolemized.  For instance, let $F_0$ be pseudo-finite, $\ch(F_0)=0$,
and let $K=F_0((t^{\qq}))$ be the field of Puiseux series over
$F_0$.   Then $K$ has Galois group $\hat{\Zz}$, and embeds into a
pseudo-finite field $F$ such that $\Aut(F^a /F) \to \Aut(K^a / K)$
is an isomorphism; hence $K$ is relatively algebraically closed in
$F$. But being Henselian and not separably closed it cannot be PAC,
by Corollary 11.5.6 of \cite{jardenfried}.
\end{Remark}

\begin{example}  A simple theory $T$ geometrically representing  $1$ and $\Zz/2\Zz$, but no other finite group.  \rm  Let $L = \{R , f,p \}$
where $R$ is a binary predicate, $p$  a unary function , $f$ a binary function symbol.   Let   $X$ denote the image of $p$,  $Y(x):= p \inv(x) \setminus \{x\}$, and $Y=\union_{x \in X} Y(x)$.  The universal theory $T_{\forall}$  then states that $R$ is a tournament on $X$,  $p: Y \to X$ has fibers of size $\leq 2$,
and $f(a,b)$ chooses an element of $Y(b)$, provided that $R(a,b)$ holds.     (Formally:    $p(p(x))=p(x)$;   $R(x,y)$ implies $px=x,py=y$, and if $px=x,py=y$ then $R(x,y)$ or $R(y,x)$ but not both;   $Y(x):= p \inv(x) \setminus \{x\}$ has at most two elements;  $f(x,y)=f(px,py)$; $pf(x,y)=y$, and $f(x,y)=y \leftrightarrow p(y)=y$.  )

It is easy to see that 
the finite $T_\forall$-structures form an amalgamation class with the joint embedding property,
and hence the model completion is a complete,
$\aleph_0$-categorical theory $T$.    Moreover, any type $tp(c/A)$ with $c \in X$ and $X(A) \neq \emptyset $
 admits an automorphism invariant extension to a universal domain $N \geq A$.  If $c \in A$ this is obvious.  If $c \notin A$ then in fact, if $d \in Y(c)$ then $tp(cd /A)$ admits   an   invariant extension to $N$.  Namely,  set $R(b,c)$ for all $b \in X(N) \setminus A$; and let
$f(b,c) =d$.
 
   Let $M \models T$, and let $a$ belong to some elementary extension $N$, with $\neg R(m,a)$ for $m \in M$.   Then the two elements of $Y(a)$    have the same type over $Ma$.   So $\Zz/2\Zz$ is geometrically represented in $T$.  
   
   We wish to show that no other groups may be geometrically represented in $T$, or even in $T^{eq}$.  For the latter,  we need to understand    algebraic closure in (potentially)  imaginary sorts.  
 Let $A=\dcl(A)$ be a finite structure, with $X(A) \neq \emptyset$.  We claim that $\acl(A) \meet (X \union Y) = \acl^{eq}(A)$, i.e. the algebraic closure
 of $A$ in the sorts $X,Y$ accounts for the algebraic imaginaries.  To see this, we may add to $A$ the elements of $p \inv(X(A))$, so
 that  $p: Y(A) \to X(A)$ is a 2-1 map.  Let $e \in  \acl(A)$ be an imaginary element.  Let $B$ be a finite set of   elements of $X \union Y$, with $e \in \dcl(A)$;
say $X(B) \setminus A =\{b_1,\ldots,b_n\}$.  Let $b_{0}$ be an element such that $R(b_0,b_i)$ holds for   $1 \leq i \leq n$.   Note that (thanks to $f$),
$Y(B) \subset \dcl(X(B) \union \{b_0\})$.   So $e \in A(b_0,b_1,\ldots,b_n)$.  But $tp(b_i/ A(b_0,\ldots,b_{i-1})$ extends to an
$A(b_0,\ldots,b_{i-1})$-invariant type over the universal domain.  Hence by reverse induction on $i$ we have $e \in A(b_0,\ldots,b_{i-1})$ for each  $i$, and so for $i=0$ we have $e \in \dcl(A)$.

On the other hand, if $M \subset A \subset N^{eq}$ then $Aut(\acl(A)/A)$ can have no more than two elements.   To see this,  we may take $A$ 
finite.  By the remark above, it suffices to consider $G=Aut( \union_{a \in X(A)} Y(a) / A)$; this group will only grow if
$A$ is restricted to the part in $X$; so say $A = \{a_1,\ldots,a_n \}$ with $a_i \in X$.  If for any $a_i$, for some $a_j$ we have $R(a_j,a_i)$,
then $G$ is trivial, because of $f$.  Otherwise, for some $a_i$, for all $a_j$ with $j \neq i$ we have $R(a_i,a_j)$.  In this case all $Y(a_j)$ with 
$j \neq i$ we have $Y(a_j) \subset \dcl(A)$, so $G$ acts faithfully on $Y(a_i)$, and hence has at most two elements.   
\end{example}  

We could easily modify this example so as to find a theory representing all subquotients of some fixed finite group $H$, but no other
finite group.

\begin{Remark}  \label{product-3}     
If two finite groups $G,H$ are geometrically represented in the complete, stable theory $T$, then so is $G \times H$. { \rm Using Remark \ref{product1}, we may take $G,H$
to be represented over the same base model $M_0$.  Say $G=\Aut(B/A)$,
$M_0 \leq A \leq B \leq N$, with $B \subset \acl(A)$ and $B$ normal over $A$;
and similarly  $H=\Aut(D/C)$.  
 It follows that $Aut(BD/AC) = G \times H$.    }  This property is inherited by the model completion of the theory  of models of $T$ with a distinguished
 automorphism, if it exists; hence it holds for $ACFA$ and for the theory of pseudo-finite fields.  
    \end{Remark}

 \begin{question}  Which finite groups can be geometrically represented in theories of pseudo-finite fields?  \end{question}

By Remark \ref{product-3} and    Theorem \ref{aut-conv} that any finite {\em Abelian} group can be geometrically represented in the theory of pseudo-finite fields containing the roots of unity.   Perhaps the internal Galois groups are indeed all Abelian.

We end with some open questions.    To state them algebraically, recall the standard description of the basic structure of the cyclotomic extension of $k=\mathbb{Q}$.    $\Aut(k(\mu_{p^{\infty}})/k)$ is the
inverse limit of the automorphism groups of the finite extensions
$\Aut(k(\mu_{p^n})/k)$ and if $p\neq 2$,
$$\Aut(k(\mu_{p^n})/k)\simeq \zz/p^{n-1}\zz\times \zz/(p-1)\zz.$$

For $i\geq j$, the restriction homomorphism

$$\begin{array}{cccc}
r_{ij}:&\Aut(k(\mu_{p^i})/k)&\longrightarrow & \Aut(k(\mu_{p^j})/k)\\
&\phi&\longmapsto & \phi_{\mid k(\mu_{p^j})},\\
\end{array}$$
which is certainly onto, respects the decomposition. Hence
$$\Aut(k(\mu_{p^{\infty}})/k) \simeq \zz_p\times \zz/q\zz.$$
where $q=p-1$;
this is also valid for $p=2$, taking $q=2$.

 Let $L_p$  be the subfield of $k(\mu_{p^{\infty}})$ fixed by $$\zz/q\zz<\zz_p\times \zz/q\zz$$
  and let $\w$ be a primitive $p$-th root of unity if $p\neq 2$ and $\sqrt{-1}$ if $p=2$.
 The field $L_p$ does not contain any $p^n$-th roots of unity except $\pm 1$. Suppose it does; then $L_p$ contains $\w$, hence
 the automorphism group of $L_p/k$ contains a subgroup of index $q$, but it is impossible since
 $\Aut(L_p/k)\simeq \zz_p$. But $L_p(\w)=k(\mu_{p^{\infty}})$ and contains $\mu_{p^{\infty}}$.

On the other hand, the finite extension $L_p(\w)=k(\mu_{p^{\infty}})$ contains all $p^n$-th roots  of unity.

 \def\Qq{\mathbb Q}
 \bigskip

 \begin{example}  \rm
 Let $F$ be a pseudo-finite field whose algebraic points intersect $\Qq(\mu_{p^{\infty}})$
 in $L_p$.  Let $F'=F(\mu_{p^{\infty}})$.  By Theorem
 \ref{aut-conv}, $p$ can be   geometrically represented in $Th(F')$.  When $p=2$,
 by Theorem \ref{aut}, $p$ cannot be geometrically represented in $Th(F)$.  \end{example}

 We   have not settled whether this phenomenon persists for $p>2$.  We formalize the question algebraically:

 Let $F$ be an algebraically closed field (say $\Qq^a$.)  Let $G= \Aut(F(t)^a / F(t))$ be the absolute Galois group of $F(t)$.
 Let $p$ be a prime, and $\mu_p^{\infty}=\union_n \mu_{p^n}$ the group of $p^n$'th roots of $1$.  Let $1^{1/(p-1)}$ be a primitive $p-1$-root of $1$.

 Let $\si \in \Aut(F)$.  Let $G(\si)$ be the centralizer in $G$ of some lifting of $\si$ to $\Aut(F(t)^a)$;
 so $G(\si)$ is a closed subgroup of $G$, determined up to conjugacy.    Let $G(\si,p)$ be the pro-$p$ part of $G(\si)$.  We have:

 \begin{enumerate}
 \item  If $\si$ fixes $\mu_p^{\infty}$, then $G(\si,p)$ is a large pro-$p$ group.  (In particular contains $\Zz_p$.)
 \item  If $Fix(\si) [1^{1/(p-1)}] \meet \mu_p^{\infty}$ is finite,    then $G(\si,p)=1$.   \end{enumerate}

 Almost all $\si$ (for the Haar measure) fall into case (ii) for all $p$, and for them we have $G(\si)=1$.

 \begin{question}  What about the intermediate cases?  In particular let $\si(x)=x^{-1}$ for $x \in \mu$.  Is $G(\si,p)=1$?  \end{question}

  \begin{problem}  Explain a priori why (if it is indeed the general case) $G(\si)$ depends only on the action
 of $\si$ on roots of unity.
  \end{problem}

\end{document}